\documentclass[fleqn,preprint,12pt,3p]{elsarticle}
\usepackage{amsmath,amssymb,amsthm,array}
\usepackage{graphicx}
\usepackage{epstopdf}
\usepackage{subcaption}
\usepackage[normalem]{ulem}
\usepackage{mathabx}

\begin{document}

\begin{frontmatter}
\title{On compactly supported discrete radial wavelets in $L^2(\mathbb{R}^2)$ and application in Tomography}
\author{K. Z. Najiya, Akshaya Ravichandran and C. S. Sastry \\
Department of Mathematics \\
Indian Institute of Technology, Hyderabad, 502285, India. \\ Email:\{ma17resch11004,ma18mscst11001 and csastry\}@iith.ac.in}
\date{}
\begin{abstract}
Radially symmetric wavelets possessing multiresolution framework are found to be useful in different fields like Pattern recognition, Computed Tomography (CT) etc. The compactly supported wavelets are known to be useful for localized operations in applications such as reconstruction, enhancement etc. In this work we introduce a novel way of designing compactly supported radial wavelets in $L^2(\mathbb{R}^2)$ from a 1D Daubechies wavelets and obtain a reconstruction formula possessing multiresolution framework. Further, we demonstrate the usefulness of our radial wavelets in Tomography. 
\end{abstract} 

\end{frontmatter}

\section{Introduction}
  
 Wavelets provide bases for several function spaces \cite{ID}. A framework [6] through which compactly supported, orthogonal (biorthogonal),  regular and real wavelets are constructed is called  multiresolution analysis (MRA). A real valued function $\eta$ defined on $\mathbb{R}^n$ is radially symmetric if $\eta (x) = \eta_1 (\| x \|),\; \forall x \in \mathbb{R}^n$ for some univariable function $\eta_1$. The wavelets, that are  radially symmetric, compactly supported and associated with multiresolution framework,  have promising applications in fields like: 
\begin{itemize}
     \item Computed Tomography \cite{radialwavelet1}\cite{madych}, where the angle-dependent wavelet weighted ramp filter can be made free from angle. As a result, direct and simplified reconstruction algorithms can be designed for different scanning geometries, and  interior tomography.
     \item Data analysis and classification \cite{sastry1}, where rotation invariant feature extraction can be realized without any approximations.
\end{itemize}
In the seminal work of \cite{madych}, it was shown that the radial wavelets possessing multiresolution framework are useful for reconstruction in Computed Tomography (CT). The author of \cite{madych} used radial functions that provide continuous wavelet transform by discretizing the scale parameter. 
 Though the work proposed the local reconstruction methods for low pass and high pass versions of function $f$ via the $\Lambda$ operator, the filter functions used for these high and low pass versions were not compactly supported. The compactly supported filters are known to be useful for localized operations in applications such as reconstruction, enhancement etc.
It can be shown (discussed in later sections) further that the compactly supported radial wavelets based on the Daubechies wavelets possess many attractive properties for numerical computations and hence are suitable for interior reconstruction in CT.  But such wavelets, to our knowledge, are not available even to date \cite{radialwavelet1}\cite{radialwavelet3}. Motivated both by the novel contributions of the ideas developed in \cite{madych} and the apparent scope for improvement through finitely supported discrete radial wavelets, we design radial wavelets possessing desired properties and show their usefulness in CT. 
   
\par The paper is organized as several  sections. In sections 2 and 3, we discuss the motivation for our work and the relevant basics of wavelets respectively. In sections 4, 5 and 6, we present  our new way of designing radial wavelets, a reconstruction as applicable in CT and simulation results  respectively. 

\section{Motivation for our work and Contribution}
The Radon transformation of a well-behaved two variable function $f$ is defined by
\begin{equation}
    R_\theta f (t)={\displaystyle\int_\mathbb{R} f(t u_\theta+s v_\theta)d s}, 
\end{equation}
which is a line integral of $f$ along the line which is perpendicular to $u_\theta$. Throughout this paper, we use $u_\theta = (\cos\theta,
\sin\theta)$ and $v_\theta = u_{\theta+\frac{\pi}{2}}$, the unit
vectors along $\theta$ and its perpendicular directions
respectively. Fourier-Slice theorem \cite{kuchment} implies that $\hat{f}(t u_\theta)=\widehat{R_{\theta}f}(t)$. Here, the ``hat" notation refers to the Fourier transform operation. A well-behaved function can be reconstructed from its Radon projections via the following inversion formula
\begin{equation*}
    f(x)={\displaystyle\int_0^{\pi}\int_{\mathbb{R}}\widehat{R_{\theta}f}(t)|t| e^{2\pi it\langle x,u_\theta\rangle}  dt d\theta},\\
\end{equation*}
 which is known as the back projection formula. The operator, $\Lambda$ defined \cite{madych} by
 $\Lambda f(x) = -\frac{1}{4\pi} \Delta R^{\#}Rf(x)$,
and its inversion $\Lambda^{-1}f(x)$ are respectively regarded as  high and low band-pass filtered versions of $f$. In practical applications of local tomography, however, one does not have to compute $\Lambda f(x)$, but attempts to reconstruct $\Lambda(\Phi_{\epsilon} \ast f)$ for some approximation of identity $\Phi_{\epsilon}(x)$ = $\epsilon^{-2} \Phi(x/\epsilon)$, where $\displaystyle\int_{\mathbb{R}^2}\Phi(x) dx=1$. It is worth a mention here that even if $\Phi$ is compactly supported, in general, $\Psi =\Lambda \Phi $ is not. In \cite{madych}, for a function that is compactly supported on unit disk, the author provided a reconstruction formula 
\begin{equation} \label{Madych_recon}
f(x) =  \frac{1}{2 \pi^2 2^{N_0} \epsilon} \int_0^{2 \pi} R_\theta f (\langle x, u_\theta \rangle ) d \theta + \sum_{j=-\infty}^{N_0} \Psi_{2^j \epsilon} \star f(x),  
\end{equation}
where, $\Psi_{2^j \epsilon} \star f(x)$ was  written in terms of $R_\theta f$ through the radial wavelet
\begin{equation}
    \Psi(x) = \frac{1}{\pi^2} \biggl( \frac{1}{\| x \|} - \frac{\sqrt{2\| x\|^2-1}}{\|x\|^2} \chi(2 \| x \|) +  \frac{\sqrt{\| x\|^2-1}}{\|x\|^2} \chi( \| x \|) \biggl). 
\end{equation}
In the above equation, $\chi$ is the characteristic function defined over $(-\infty,-1]\cup[1,\infty)$. Consequently, $\Psi$ is not compactly supported. As a result, the numerical  computation of  $\Psi_{2^j \epsilon} \star f(x)$ at every $x$ requires all values of $f(x)$. In contrast, a finitely supported radial $\Psi$ would only require the values of $f$ from a small disc around $x$ of radius dependent on the radius of support of $\Psi_{2^j \epsilon}$. Consequently, a reconstruction formula of the type (\ref{Madych_recon}) in terms of Daubechies discrete wavelets looks more promising, as such wavelets possess many attractive properties such as finite support, vanishing moments, enough smoothness etc, which are very useful for numerical computations and local recovery. Motivated by the stated points, we intend to define compactly supported radially symmetric wavelets and obtain a discrete wavelet based multiresolution formula that is capable of being used for recovering $f$ precisely in such situations as local recovery.  Though the local recovery was addressed through conventional wavelet decomposition (wherein each coefficient is computed via the backprojection algorithm) and other frameworks   \cite{sas_das}\cite{radialwavelet1}\cite{abcd}\cite{klan} (and the references therein), the proposed radial based approach results in a simpler formula involving angle-independent ramp filter.  While the basic philosophy of the present work is inspired by the ideas in \cite{madych}, our contribution lies in designing and using novel compactly supported radial wavelets obtained through the known 1D wavelets that possess attractive numerical properties. 

\section{Basics of Wavelets}

\par A discrete wavelet in $L^2(\mathbb{R})$ is a function $\psi\in  L^2(\mathbb{R})$ such that the collection $\{\psi_{j,m}\}_{j,m\in \mathbb{Z}}$ forms an orthonormal basis for $L^2(\mathbb{R})$, where $\psi_{j,m}(t)=2^{j/2}\psi(2^j t-m)$. Multi Resolution Analysis (MRA) is a framework through which several classes of wavelets can be constructed. 
The approximation and details spaces are defined as 
$V_j=\overline{span}\{ \phi_{j,k}: k \in \mathbb{Z} \}$, 
$W_j= V_{j+1} \ominus V_j = \overline{span}\{ \psi_{j,k}: k \in \mathbb{Z} \}$ 
where $\phi$ and $\psi$ are the scaling and wavelet functions associated with the Daubechies class of wavelets \cite{ID}. The 2D wavelets generated via a tensor product of 1D wavelets are found useful in many applications. 
As highlighted in the previous section, radial wavelets do have a potential for local reconstruction in tomography. Despite this, to the best of our knowledge, the compactly supported real radial wavelets associated with the multiresolution framework are not known to exist.  Nevertheless, a 2D  function $\eta^{(2)}$, defined from a 1D compactly supported locally integrable function $\eta^{(1)}$ as
\begin{equation} \label{rad}
\eta^{(2)}(x) = \frac{1}{2 \pi} \int_0^{ 2 \pi} \eta^{(1)}(<x, u_\theta>)\eta^{(1)}(<x, v_\theta>) d \theta\\
\end{equation}    
is radially symmetric as $\eta^{(2)}(A_\theta x)=\eta^{(2)}(x)$ for any rotation matrix $A_\theta$. It may be noted that $\eta^{(2)}$ is compactly supported 
in the disk centered at origin with radius $\sqrt{2} a$ when the support of $\eta^{(1)}$ is $[-a, a]$.  
In Fourier domain, however, $\eta^{(2)}$ satisfies
$$\widehat{\eta^{(2)}}(\omega)= \frac{1}{2\pi} {\displaystyle\int_0^{2\pi} \widehat{\eta^{(1)}}(\langle \omega,u_\theta \rangle) \widehat{\eta^{(1)}}(\langle \omega,v_\theta \rangle) d\theta}.$$
\noindent Making use of the stated definition, we define compactly supported 2D wavelets and show their usefulness towards reconstruction in CT.
\section{A reconstruction formula via radial wavelets in $L^2(\mathbb{R}^2)$} \label{sec:Radial_Recon}
Starting with a compactly supported, sufficiently smooth, real scaling and wavelet pair  $\phi,\psi$ (associated with an MRA framework) in $L^2(\mathbb{R})$, we define  the auto correlation function of $\phi$ and $\psi$, respectively,  as $\Phi=\phi* \tilde{\phi}$ and $\Psi=\psi* \tilde{\psi}$ where $\tilde{\phi}(x)=\phi(-x)$. Now the radial scaling and radial wavelet functions are defined in auto-correlation domain as
\begin{equation}
    \Phi^r(x) = \frac{1}{2\pi} {\displaystyle\int_0^{2\pi} \Phi(\langle x,u_\theta \rangle)\Phi(\langle x,v_\theta \rangle) d\theta} \text{ and } 
    \Psi^r(x) = \frac{1}{2\pi} {\displaystyle\int_0^{2\pi} \Psi(\langle x,u_\theta \rangle)\Psi(\langle x,v_\theta \rangle) d\theta}.
\end{equation}
It follows that the compactly supported and radially symmetric $\Phi^r(x)$ and $\Psi^r(x)$ satisfy
\begin{equation}\label{eq:mean_radial}
\begin{split}
\int_{\mathbb{R}^2} \Phi^r(x) dx  =1, 
\int_{\mathbb{R}^2} \Psi^r(x) dx  = 0. 
\end{split}
\end{equation}
By 1-D MRA framework \cite{ID}, we have
$|\widehat{\phi}(\omega)|^2+\sum_{j=0}^{\infty}{ |\widehat{\psi}(2^{-j}\omega)|^2}=1 \; a.e \;  \omega \in \mathbb{R}.$
With the notation:  
$\phi_j(x)= 2^j \phi(2^jx)$, we have $|\widehat{\phi}(2^j\omega)|^2= {\widehat{\Phi}}_{-j}(\omega)= {\widehat{\Phi}}(2^j\omega)$. Consequently, for any $J \in \mathbb{Z}$, we have
\begin{equation}\label{eq:form1}
{\widehat{\Phi}}_J(\omega)+ \sum_{j=J}^{\infty}{ \widehat{{\Psi}}_j(\omega)}=1 \; a.e \; \omega \in \mathbb{R}.
\end{equation}
Substituting $\langle \xi,u_\theta \rangle$,  $\langle \xi,v_\theta \rangle$ for $\omega$ in (\ref{eq:form1}) separately 
and multiplying the equations, we get
\begin{equation} \label{eq:formula_midpt}
\begin{split}
{\widehat{\Phi}}_J(\langle \xi,u_\theta \rangle){\widehat{\Phi}}_J(\langle \xi,v_\theta \rangle) &+ \sum_{j=J}^{\infty}{{\widehat{\Phi}}_J(\langle \xi,u_\theta \rangle) \widehat{{\Psi}}_j(\langle \xi,v_\theta \rangle)} + \sum_{j=J}^{\infty}{{\widehat{\Psi}}_j(\langle \xi,u_\theta \rangle) \widehat{{\Phi}}_J(\langle \xi,v_\theta \rangle)} \\ &+\sum_{j=J}^{\infty}\sum_{k=J}^{\infty}{\widehat{{\Psi}}_j(\langle \xi,u_\theta \rangle) \widehat{{\Psi}}_k(\langle \xi,v_\theta \rangle)}=1; \; \;  \xi \in \mathbb{R}^2.
\end{split}
\end{equation}
Finally, integrating both sides of (\ref{eq:formula_midpt}) from 0 to $2 \pi$, multiplying with $\widehat{f}$ on both sides and then taking the inverse Fourier transform, we get
\begin{equation}\label{decomp_formula}
    f= f \star \Phi_{J,J}^{r}+2 \sum_{j=J}^\infty {f \star \Psi_{J,j}^{r,h}}+\sum_{j,k=J}^\infty {f \star \Psi_{j,k}^{r,d}},
\end{equation}
where, the radial functions $\Phi^r$, $\Psi^{r,h}$, $\Psi^{r,v}$ and $\Psi^{r,d}$ are defined as follows:
 \begin{equation} \label{4radial_fun}
\begin{split}
\Phi_{j,j}^{r}(x) &=\frac{1}{2\pi} {\displaystyle\int_0^{2\pi} {\Phi}_j(\langle x,u_\theta \rangle) {\Phi}_j(\langle x,v_\theta \rangle) d\theta} \\
\Psi_{j,k}^{r,h}(x) &=\frac{1}{2\pi} {\displaystyle\int_0^{2\pi} {\Phi}_j(\langle x,u_\theta \rangle) {\Psi}_k(\langle x,v_\theta \rangle) d\theta} \\
\Psi_{j,k}^{r,v}(x) &=\frac{1}{2\pi} {\displaystyle\int_0^{2\pi} {\Psi}_j(\langle x,u_\theta \rangle) {\Phi}_k(\langle x,v_\theta \rangle) d\theta} \\
\Psi_{j,k}^{r,d}(x) &=\frac{1}{2\pi} {\displaystyle\int_0^{2\pi} {\Psi}_j(\langle x,u_\theta \rangle) {\Psi}_k(\langle x,v_\theta \rangle) d\theta}.
\end{split}
\end{equation}
In view of the radial functions in (\ref{4radial_fun}) being defined from the auto-correlation (which are symmetric) functions of $\phi$ and $\psi$,  a simple verification implies that $\Psi_{j,k}^{r,h}(x) = \Psi_{k,j}^{r,v}(x)$.

\begin{figure}
    \centering
    \subcaptionbox{}{\includegraphics[width=2.5cm,height=2.75cm]{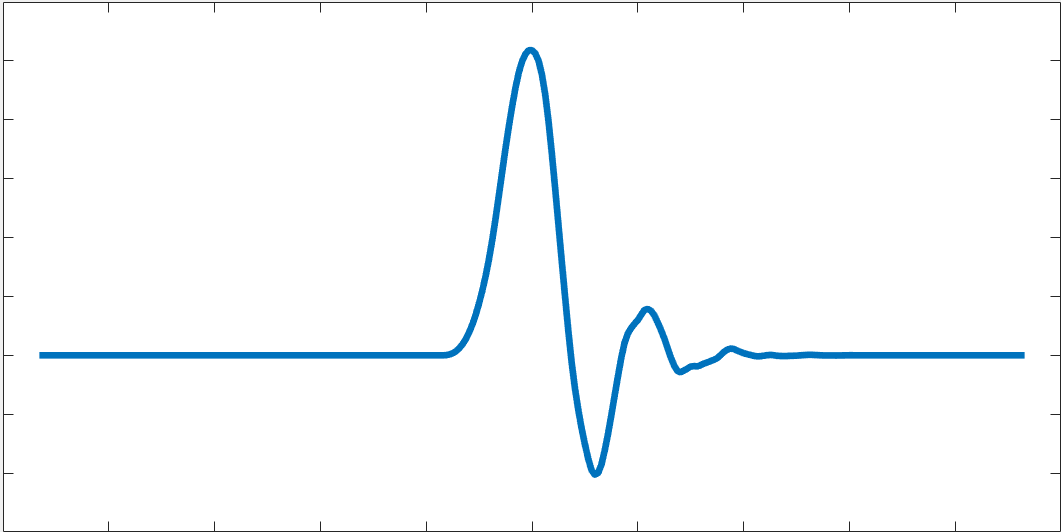}}
    \subcaptionbox{}{\includegraphics[width=2.5cm,height=2.75cm]{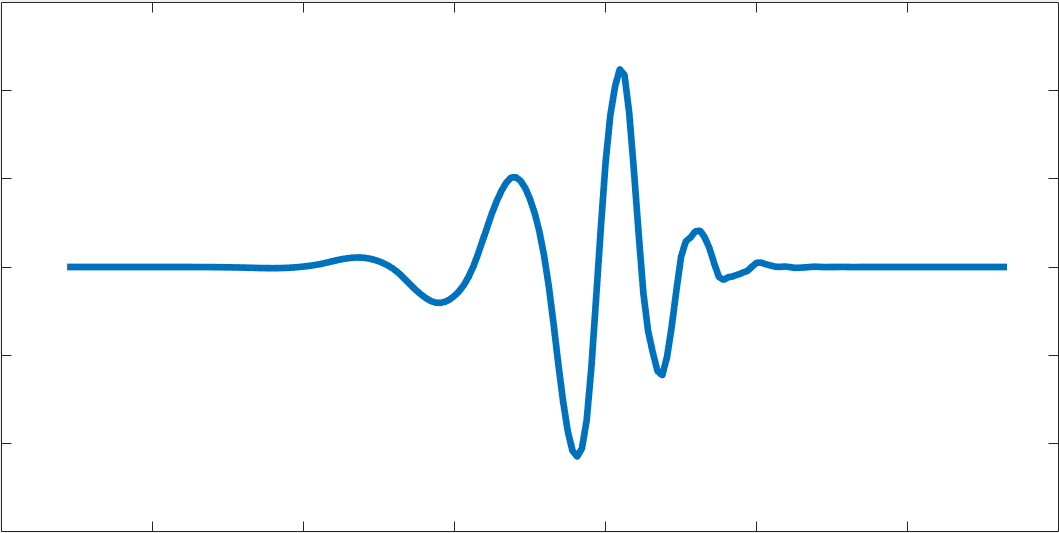}}
    \subcaptionbox{}{\includegraphics[width=2.75cm,height=2.75cm]{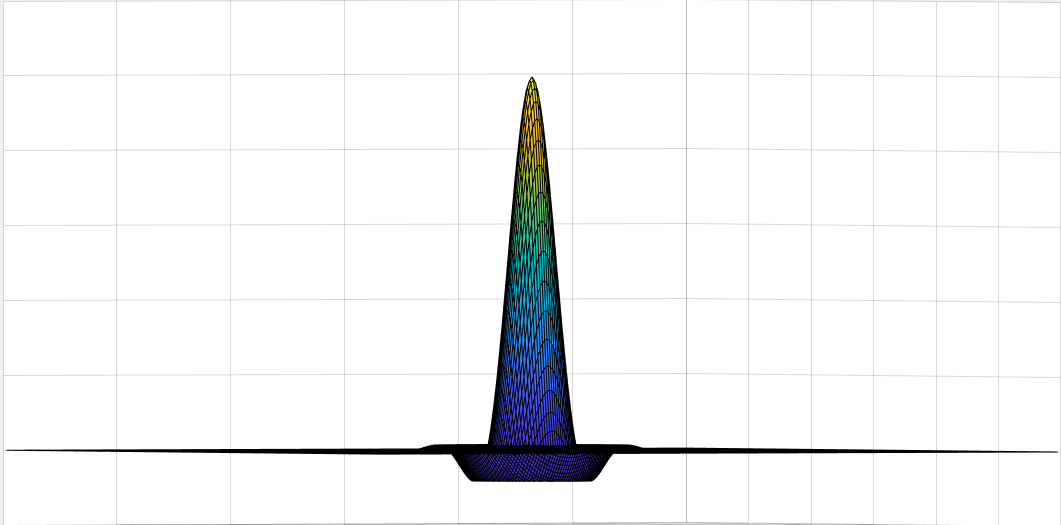}}
    \subcaptionbox{}{\includegraphics[width=2.75cm,height=2.75cm]{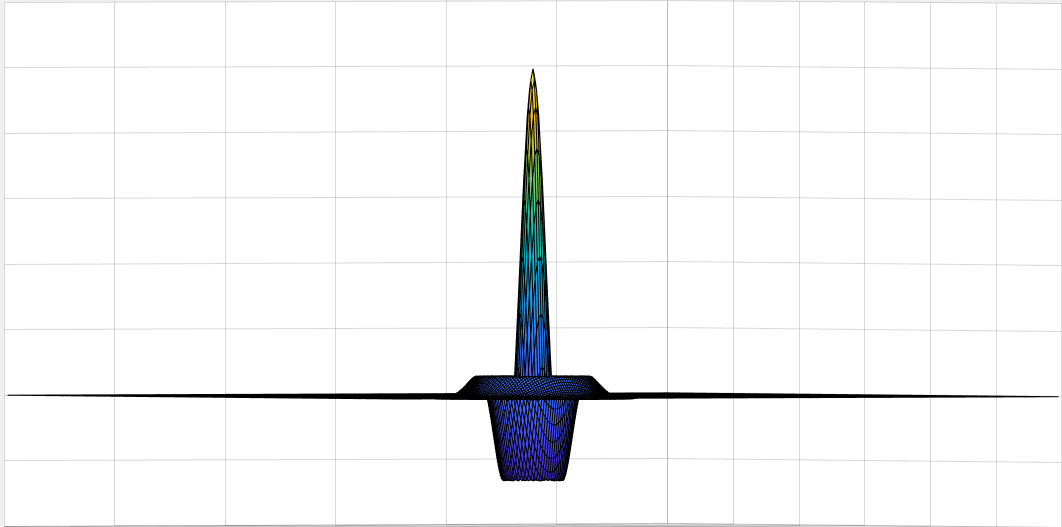}}
    \subcaptionbox{}{\includegraphics[width=2.75cm,height=2.75cm]{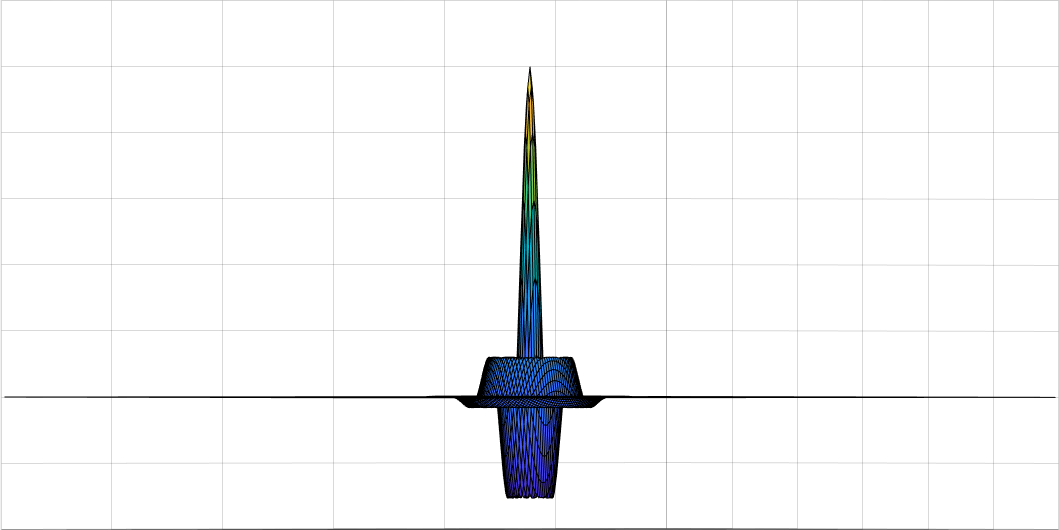}}
  \caption{(a). The 1D (a) scaling and (b).  wavelet pair associated with the Daubechies Wavelet `db6'. The 2D radial functions stated in (\ref{4radial_fun}): (b). $\Phi^r$, (c). $\Psi^{r,h}$ and (d). $\Psi^{r,d}$, which are based on (a) and (b).}
    \label{plot:radial}
\end{figure}

Like in decimal expansion of numbers, in  (\ref{decomp_formula}), details (i.e., $f\star\eta$, $\eta=\Psi_{J,j}^{r,h}, \; \Psi_{J,j}^{r,d}$) are added to the coarse approximation (i.e, $f \star \Phi_{J,J}$) to improve upon the accuracy in representation.  
The radial scaling and wavelet functions given in (\ref{4radial_fun}) are shown in Figure \ref{plot:radial}. 
\noindent
\section{CT reconstruction using radial wavelets}
In this section, we derive a formula for the computation of $f \star \eta$, where $\eta$ is any of $\Phi_{J,J}$, $\Psi_{J,j}^{r,h}$ and $\Psi_{J,j}^{r,d}.$ Consequently, the formula in (\ref{decomp_formula}) can be completely written in terms of Radon projections. For $f,\eta \in L^2(\mathbb{R}^2)$, we have
\begin{equation} \label{eq:conv_filt}
    \begin{split}
       f\star\eta(x)&={\displaystyle\int_{\mathbb{R}^2} \widehat{f\star\eta}(w) e^{2\pi ix.w}  dw}\\
       &={\displaystyle\int_0^{\pi}\int_{\mathbb{R}}\widehat{R_{\theta}f}(t)|t|\hat{\eta}(t u_0) e^{2\pi it(x. u_\theta)}  dt d\theta}
      ={\displaystyle\int_0^{\pi} [{R_{\theta}f\star (\widecheck{|.|\hat{\eta}}(. u_0)})]  (x.u_\theta) d\theta}.\\
    \end{split}
\end{equation}
In the above chain of equations, the ramp-filtered function  $\widecheck{(|.|\hat{\eta}})$ refers to the inverse Fourier transform of $(|.|\hat\eta)$, which is angle-independent due to the radial symmetry of wavelets.  Using (\ref{eq:conv_filt}), the reconstruction formula  in (\ref{decomp_formula}) can be used for CT reconstructions. The stated angle-independence makes this approach suitable even for divergent beam CT. Since our radial-wavelet based approach results in a backprojection type formula, an efficient and faster way of inversion method \cite{FastBP} can be employed. 
\par Recovering a local portion of an image from localized Radon projections in X-ray
CT has acquired importance as it results  in reduction of the radiation dosage and time of exposure. Several authors developed wavelet based methods \cite{klan} (and the references therein). But the equation (\ref{decomp_formula})
through (\ref{eq:conv_filt}) provides a simple multilevel backprojection formula. When $\hat{\eta}(0) \neq 0$, the modulus function $|.|$ introduces discontinuity in $(|.|\hat{\eta})$ thereby resulting in the spread of support of $\widecheck{(|.|\hat{\eta}})$. As a result and in view of  (\ref{eq:mean_radial}), $\widecheck{(|.|\hat{\Phi^r}})$ is no more compactly supported. Nevertheless, due to the radial structure of $\Phi^r$ coupled with the properties of associated Daubechies scaling function, the energy $\frac{\|\widecheck{(|.|\hat{\Phi^r}})\|_{L^2(\mathbb{R}^2-supp(\Phi^r))}}{\|\widecheck{(|.|\hat{\Phi^r}})\|_{L^2(\mathbb{R}^2)}} \times 100 $ is only $1.556 \times 10^{-5}$, which implies that numerically $\widecheck{(|.|\hat{\Phi^r}})$ behaves like a finitely supported function.
Consequently every term of (\ref{decomp_formula}) can be computed via 
(\ref{eq:conv_filt}) using localized data. 
\section{Simulation results}
\begin{figure}
    \centering
    \subcaptionbox{}{\includegraphics[width=2.5cm]{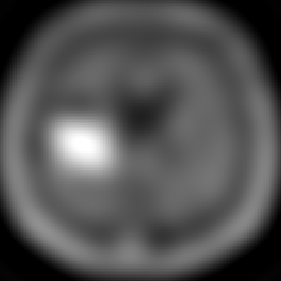}}
    \subcaptionbox{}{\includegraphics[width=2.5cm]{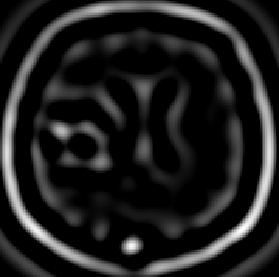}}
    \subcaptionbox{}{\includegraphics[width=2.5cm]{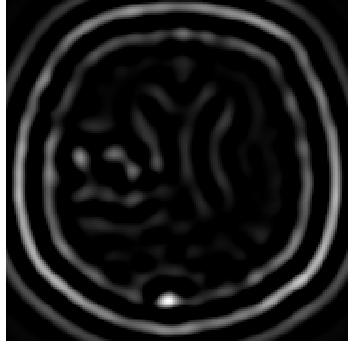}}
    \subcaptionbox{}{\includegraphics[width=2.5cm]{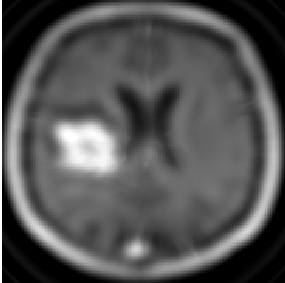}}
    \subcaptionbox{}{\includegraphics[width=2.5cm]{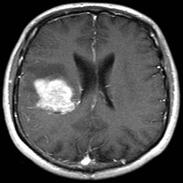}}
    \caption{For the test image in (e), the images in (a), (b) and (c) refer respectively to $f \star \Phi^r_{0,0}$, $f \star \Psi^{r,h}_{0,0}$ and $f\star\Psi^{r,d}_{0,0}$. The image in (d) is the combination of the ones in (a), (b) and (c) as in (\ref{decomp_formula}). } \label{fig:MRI_image_recon}
\end{figure}

\begin{figure}
    \centering
    \subcaptionbox{}{\includegraphics[width=2.5cm]{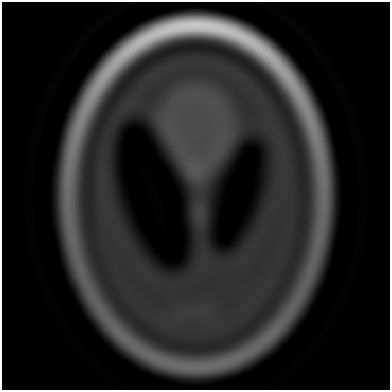}}
    \subcaptionbox{}{\includegraphics[width=2.5cm]{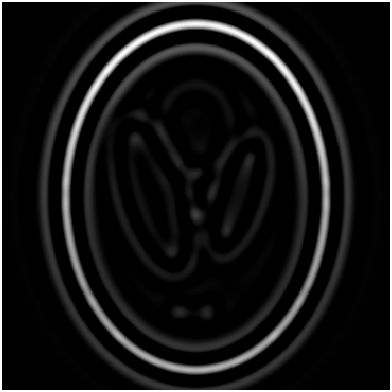}}
    \subcaptionbox{}{\includegraphics[width=2.5cm]{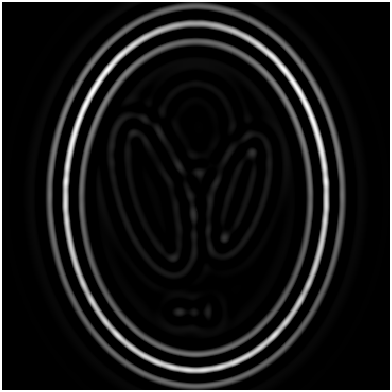}}
    \subcaptionbox{}{\includegraphics[width=2.5cm]{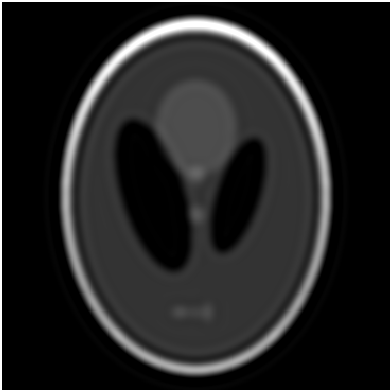}}
    \subcaptionbox{}{\includegraphics[width=2.5cm]{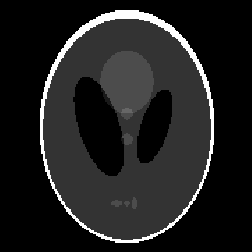}}
    \caption{For the test image in (e), the images in (a), (b) and (c) refer respectively to $f \star \Phi^r_{0,0}$, $f \star \Psi^{r,h}_{0,0}$ and $f\star\Psi^{r,d}_{0,0}$. The image in (d) is  the combination of the ones in (a), (b) and (c) as in (\ref{decomp_formula}).  } \label{fig:Phant_image_recon}
\end{figure}

We carried out our simulation work on a test image (Figure \ref{fig:MRI_image_recon}(f)) of size $183 \times 183$. The low pass and high pass images, at scale $J=0$,  reconstructed from the Radon projections, are shown in Figure \ref{fig:MRI_image_recon}. The relative 2-norm error (that is, $\frac{\|f-f_r \|_2}{\| f \|_2}$, where $f$ and $f_r$ are the original and reconstructed images respectively) in reconstruction is 0.1285. It may be noted that the reconstruction quality can be improved by adding more detail images at higher levels. Similar reconstruction results with Phantom test image are shown in Figure \ref{fig:Phant_image_recon}. 
\par
\begin{figure}
\centering
\subcaptionbox{}{\includegraphics[width=3.5cm]{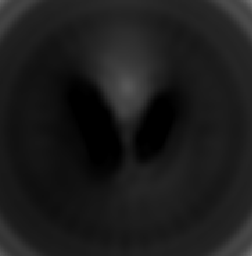}}
\subcaptionbox{}{\includegraphics[width=3.5cm]{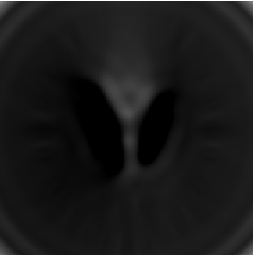}}
\caption{Local reconstruction for phantom  with  (a) $-1 = J \leq j,k \leq 1$  (b) $-1 = J \leq j,k \leq 2$ in (\ref{decomp_formula}). This figure indicates that with more detail images, the reconstruction quality increases.}
\label{plot:local}
\end{figure}

 We carried out the local reconstruction by considering a circle of radius 16 pixels around the center of the image as the interior portion. In view of the convolution structure in (\ref{decomp_formula}), which is computed via (\ref{eq:conv_filt}), we have supplied slightly excess data of 8 pixels (accounting for the width of the filter function) and obtained the local reconstruction (for the Phantom test image Figure \ref{fig:Phant_image_recon}(f))  as shown in Figure  \ref{plot:local}, which justifies that our reconstruction formula is suitable for local reconstruction from localized projections. \\

\noindent

\noindent{\bf Acknowledgments}:\\
The first author is thankful to the UGC, Govt. of India (JRF/2016/409284) for its  support. 

\normalem
\bibliography{bibliography}
\bibliographystyle{plain}

\end{document}